\documentclass{article}
\usepackage{graphicx}

\newcommand{\smalllineskip}{\baselineskip=10pt}

\def\qed{\hbox{${\vcenter{\vbox{                        
   \hrule height 0.4pt\hbox{\vrule width 0.4pt height 6pt
   \kern5pt\vrule width 0.4pt}\hrule height 0.4pt}}}$}}

\title{{\normalsize \textbf{UNKNOTTING VIRTUAL KNOTS WITH GAUSS DIAGRAM FORBIDDEN MOVES}}}
\author{{\small SAM NELSON} \\
       \textit{\small Department of Mathematics, Louisiana State University} \\
       \textit{\small Baton Rouge, LA 70803} \\
       \texttt{\small nelson@math.lsu.edu} \\
       }
\date{ }

\begin{document}
\pagestyle{empty}

\maketitle
\thispagestyle{empty}


\begin{center}
\textsc{abstract}
\end{center}
\vspace{-0.10in}

{\small \leftskip=.5in \rightskip=.5in \smalllineskip {The forbidden 
moves can be combined with Gauss diagram Reidemeister moves to obtain 
move sequences with which we may change any Gauss diagram (and hence 
any virtual knot) into any other, including in particular the unknotted 
diagram.
}

\leftskip=0in 
\rightskip=0in

\medskip
\normalsize

In 1996 Kauffman [1] introduced the theory of virtual knots, 
extending the topological concept of ``knots'' to include general 
Gauss codes. In 1999 Goussarov, Polyak and Viro [2] described 
virtual knots in terms of Gauss diagrams, which provide a visual way 
to represent Gauss codes.

Consider a classical knot diagram $K\subset \mathbf{R}^2$ as an 
immersion $K: S^1 \to \mathbf{R}^2$ of the circle in the plane with 
crossing information specified at each double point. A \textbf{Gauss 
diagram} for a classical 
knot diagram is an oriented circle considered as the preimage of the 
immersed circle with chords connecting the preimages of each double point.
We specify crossing information on each chord by directing the chord toward 
the undercrossing point and decorating each with with signs specifying the 
local writhe number.

\begin{figure}[h!]
$$
\raisebox{-0.5in}{\includegraphics[width=1in,height=1in]{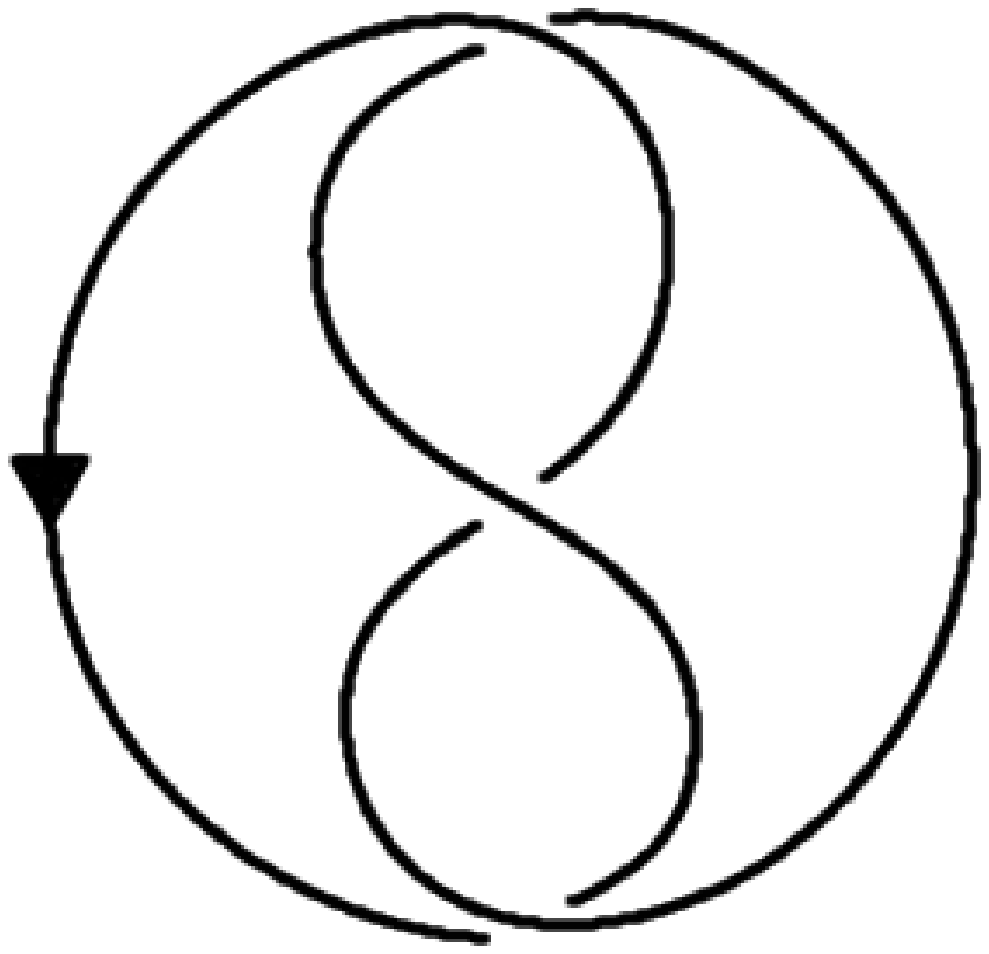}}
\raisebox{-0.5in}{\includegraphics[width=1in,height=1in]{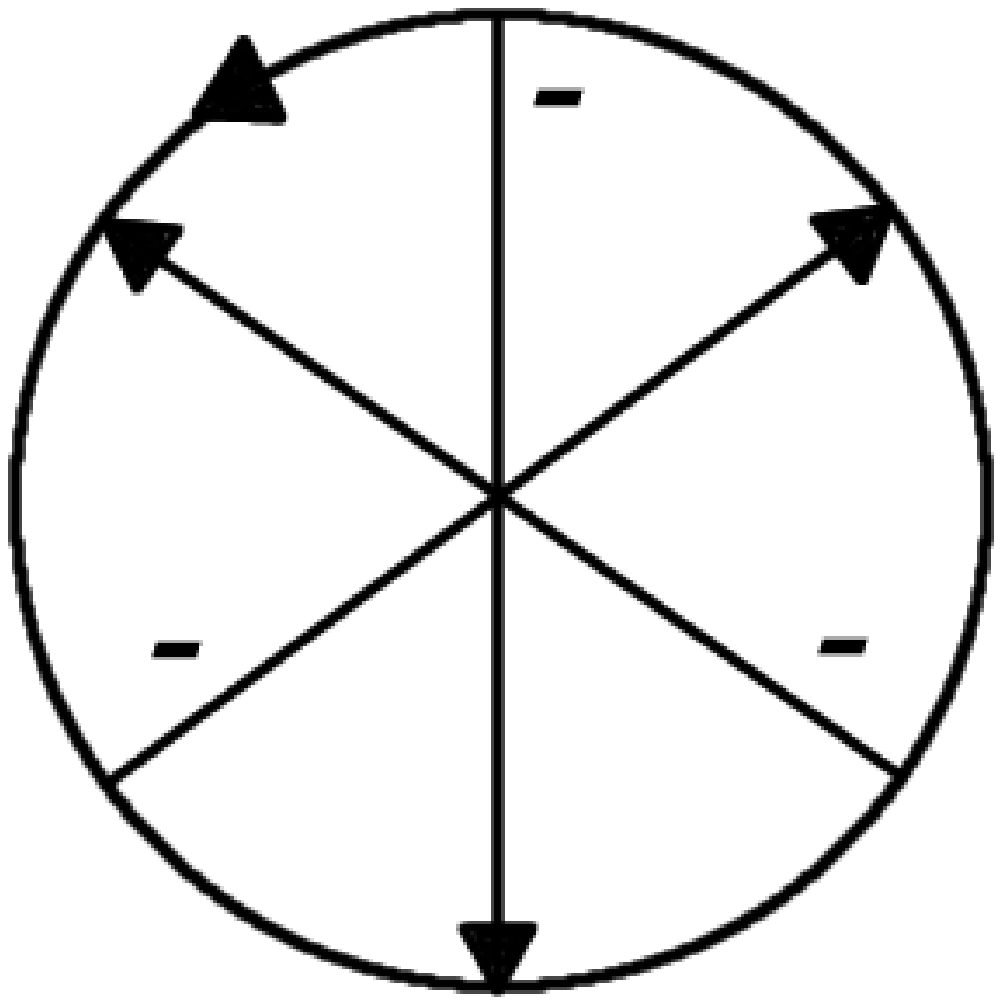}}
$$
\caption{A simple knot and its Gauss Diagram}
\end{figure}

A \textbf{virtual knot} is an equivalence class of Gauss diagrams under the
relations in Figure 2, which are the classical Reidemeister moves written in 
terms of Gauss diagrams. Note that there are several variations of move
III depending on the orientations of the strands; we only depict
two.

\begin{figure}[h!]
$$
\raisebox{-0.5in}{\includegraphics[width=1in,height=1in]{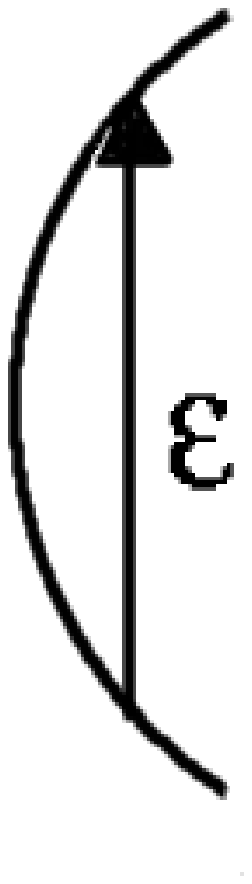}}
{\longleftrightarrow \atop \mathrm{I}}
\raisebox{-0.5in}{\includegraphics[width=1in,height=1in]{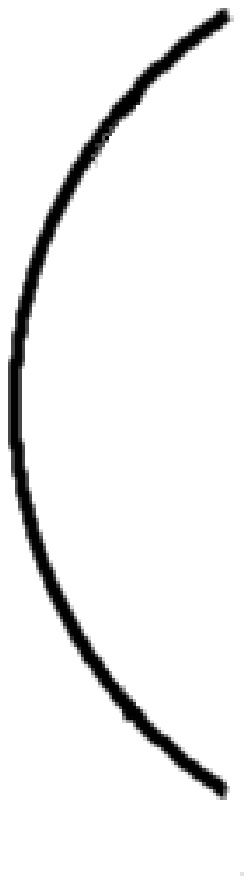}}
\raisebox{-0.5in}{\includegraphics[width=1in,height=1in]{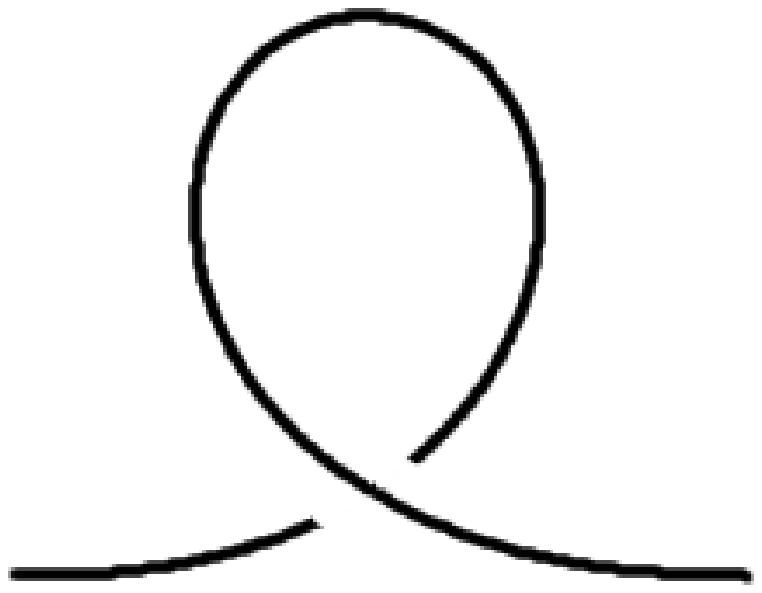}}
{\longleftrightarrow \atop \mathrm{I}}
\raisebox{-0.5in}{\includegraphics[width=1in,height=1in]{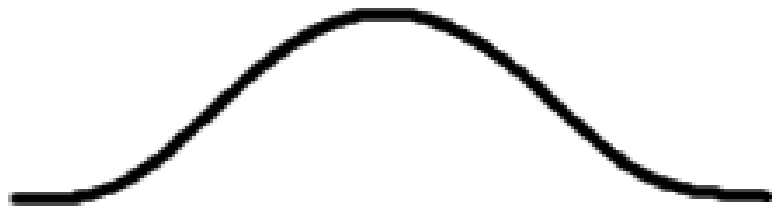}}
$$
$$
\raisebox{-0.5in}{\includegraphics[width=1in,height=1in]{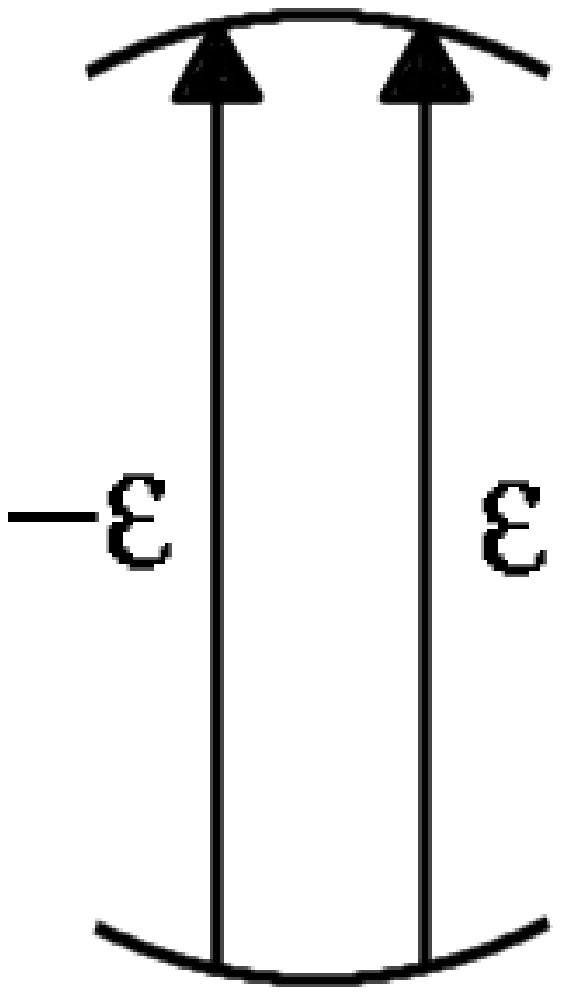}}
{\longleftrightarrow \atop \mathrm{II}}
\raisebox{-0.5in}{\includegraphics[width=1in,height=1in]{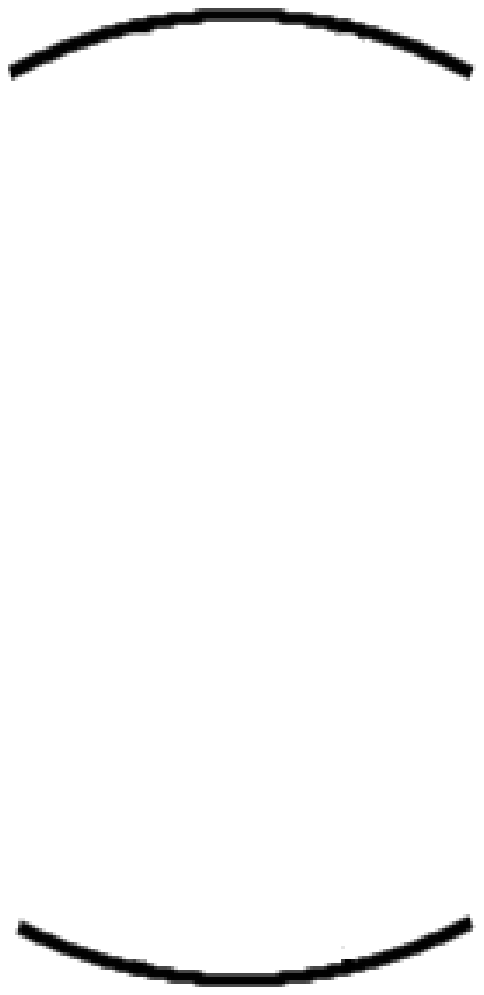}}
\raisebox{-0.5in}{\includegraphics[width=1in,height=1in]{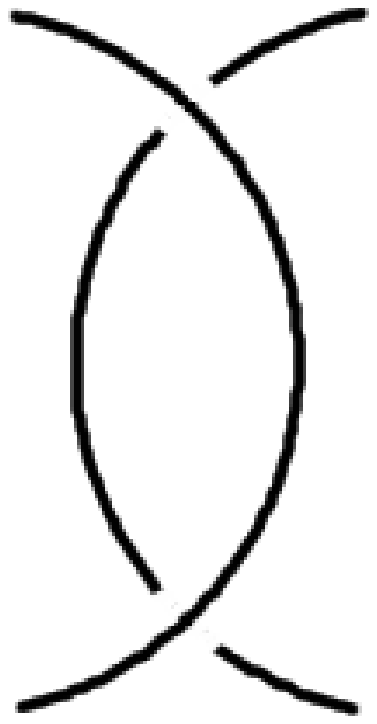}}
{\longleftrightarrow \atop \mathrm{II}}
\raisebox{-0.5in}{\includegraphics[width=1in,height=1in]{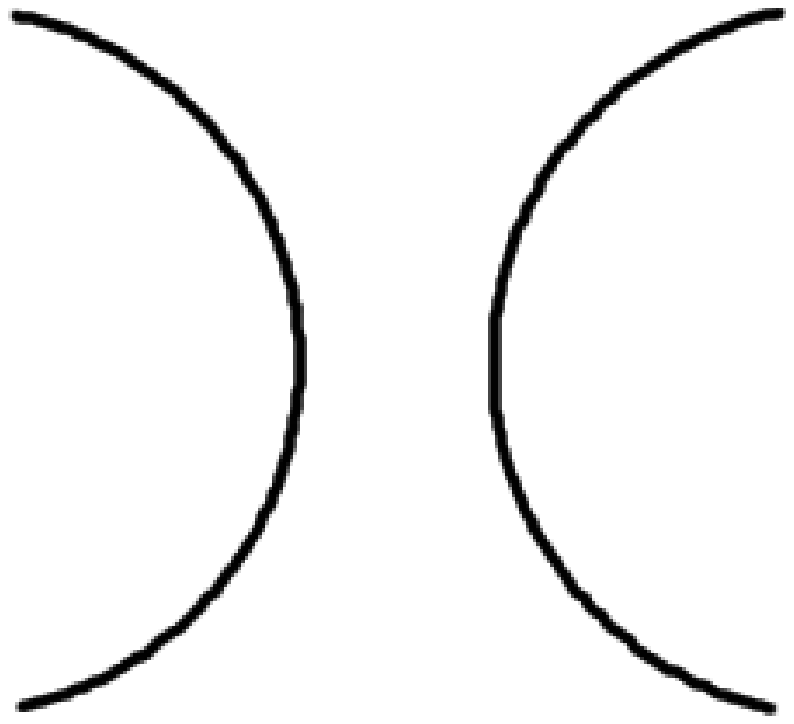}}
$$
$$
\raisebox{-0.5in}{\includegraphics[width=1in,height=1in]{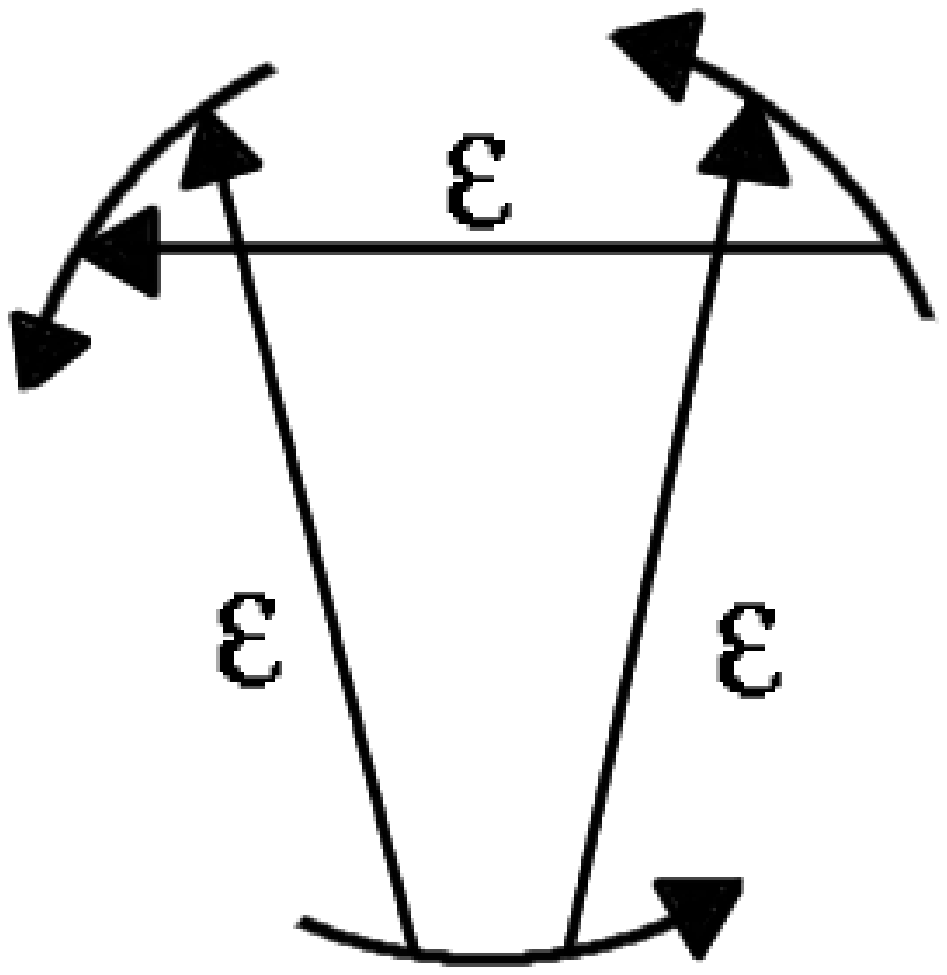}}
{\longleftrightarrow \atop \mathrm{III}}
\raisebox{-0.5in}{\includegraphics[width=1in,height=1in]{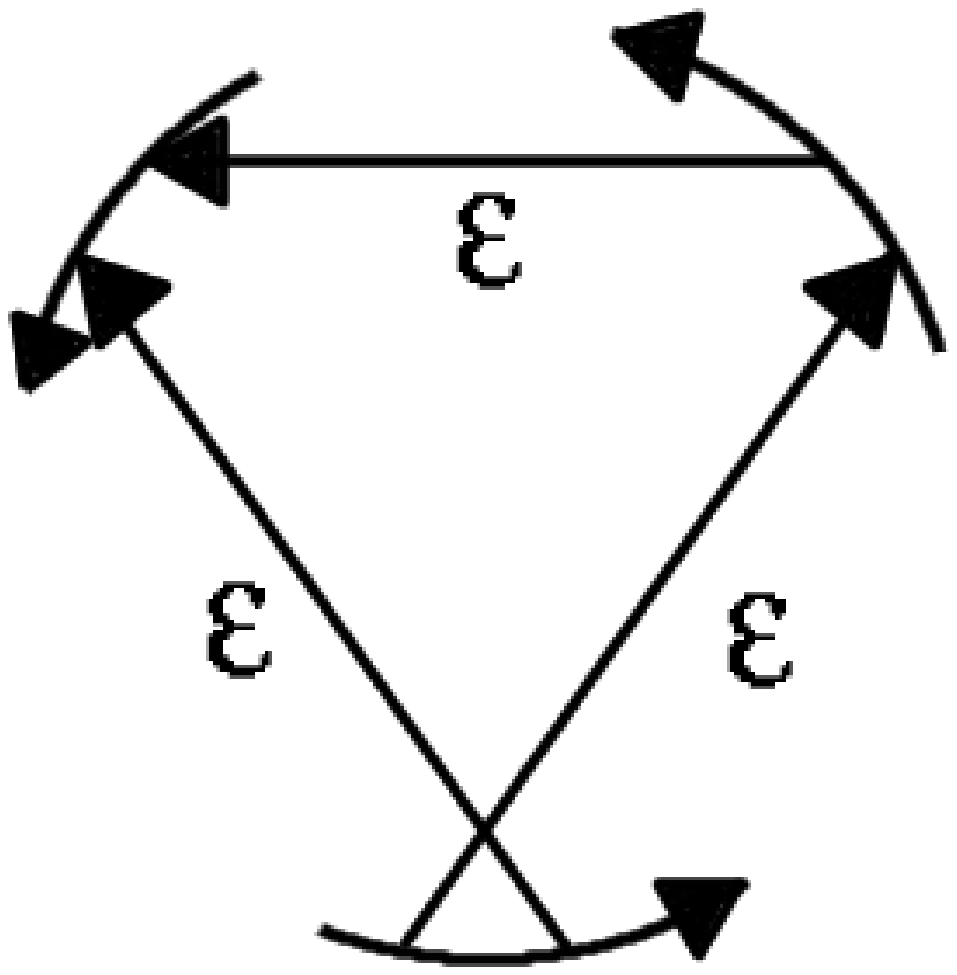}}
\raisebox{-0.5in}{\includegraphics[width=1in,height=1in]{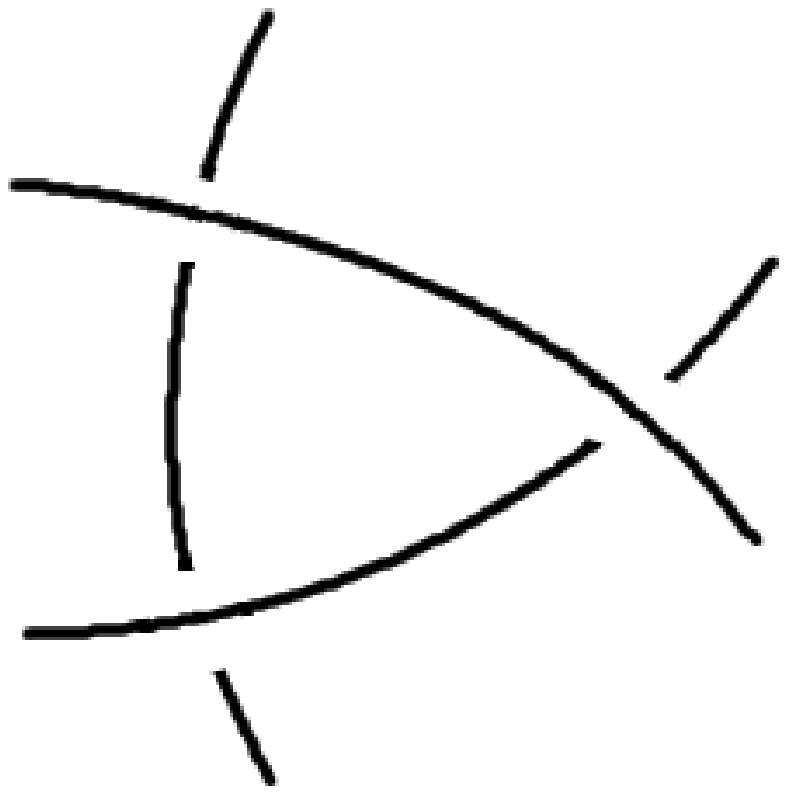}}
{\longleftrightarrow \atop \mathrm{III}}
\raisebox{-0.5in}{\includegraphics[width=1in,height=1in]{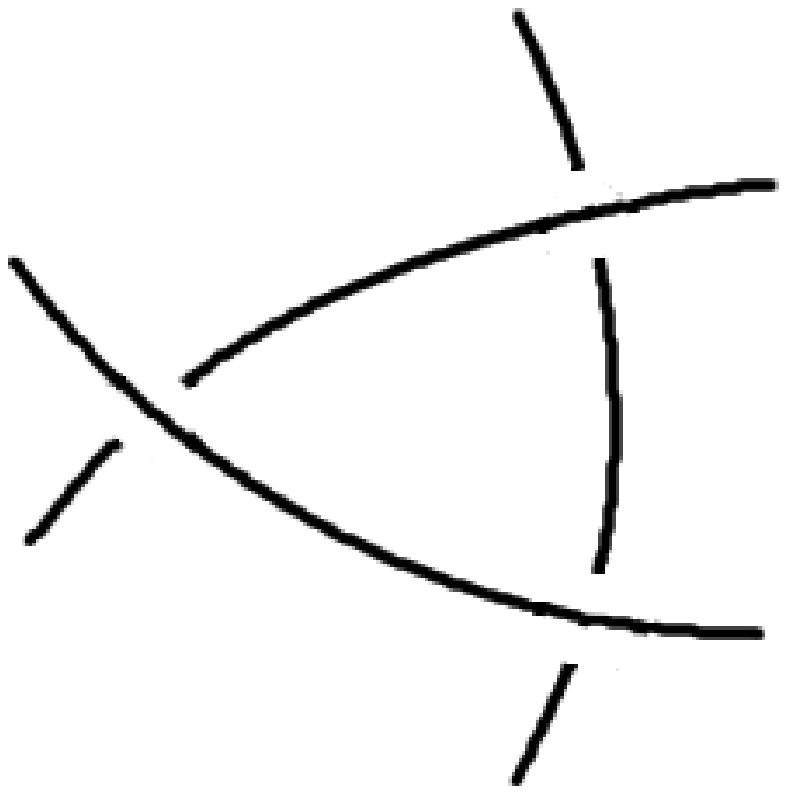}}
$$
\caption{Moves I, II and III.}
\end{figure}

Not every Gauss diagram corresponds to a classical knot. Indeed, a Gauss 
diagram determines a 4-valent graph with crossing information specified
at the vertices; such a graph represents a classical knot or link diagram if 
it is planar, but, of course, not every 4-valent graph is planar. To draw 
non-planar graphs in the plane, we usually introduce crossings, but 
these new crossings must be kept distinct from the vertices, which represent
classical crossings specified by the Gauss diagram. To draw these 
non-planar graphs as \textbf{virtual knot diagrams}, we 
introduce \textbf{virtual crossings} to distinguish crossings arising from 
non-planarity of the graph from real crossings represented by vertices. 
Virtual crossings are drawn as an intersection surrounded by a circle. A 
virtual crossing has no over- or under-sense and no sign, and virtual 
crossings do not appear on Gauss diagrams -- they are artifacts of our 
representing a non-planar graph in the plane.

A Gauss diagram determines a neighborhood of each real crossing and the 
order in which the edges entering and leaving such a neighborhood are 
connected. Outside these neighborhoods, we are free to draw the arcs 
connecting the neighborhoods however we want, introducing virtual crossings 
as necessary. The \textbf{virtual moves} in figure 3 allow us to change 
any virtual knot diagram representing a particular Gauss diagram into any 
other virtual knot diagram representing the same Gauss diagram by allowing
the interior of an arc containing only virtual crossings to be moved 
arbitrarily around the diagram.

\begin{figure}[h!]
$$
\raisebox{-0.5in}{\includegraphics[width=1in,height=1in]{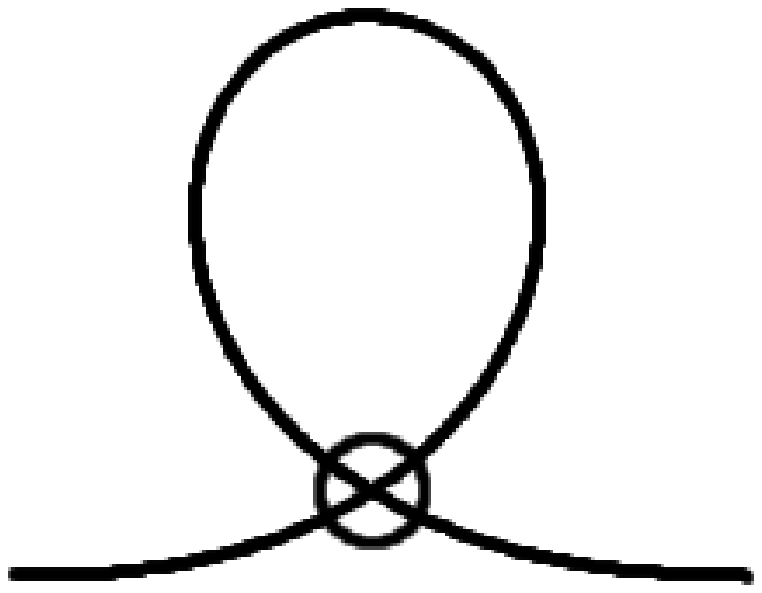}}
{\longleftrightarrow \atop v\mathrm{I}}
\raisebox{-0.5in}{\includegraphics[width=1in,height=1in]{i4.eps}}
\raisebox{-0.5in}{\includegraphics[width=1in,height=1in]{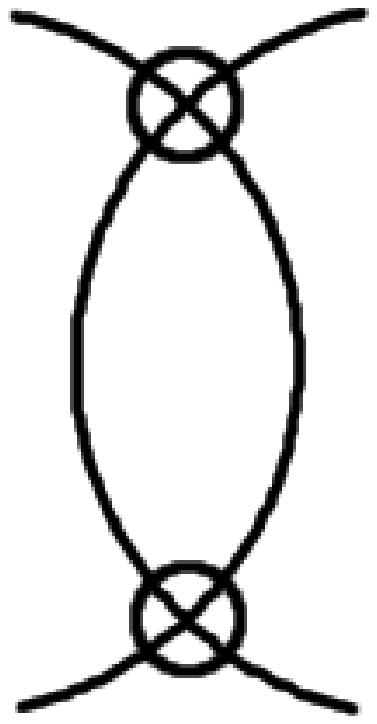}}
{\longleftrightarrow \atop v\mathrm{II}}
\raisebox{-0.5in}{\includegraphics[width=1in,height=1in]{ii4.eps}}
$$
$$
\raisebox{-0.5in}{\includegraphics[width=1in,height=1in]{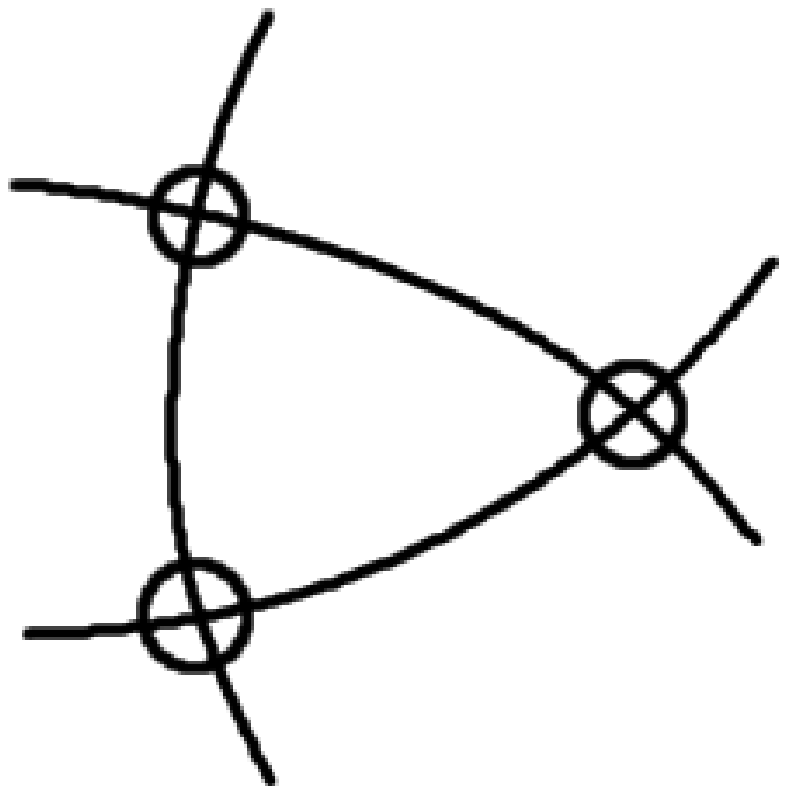}}
{\longleftrightarrow \atop v\mathrm{III}}
\raisebox{-0.5in}{\includegraphics[width=1in,height=1in]{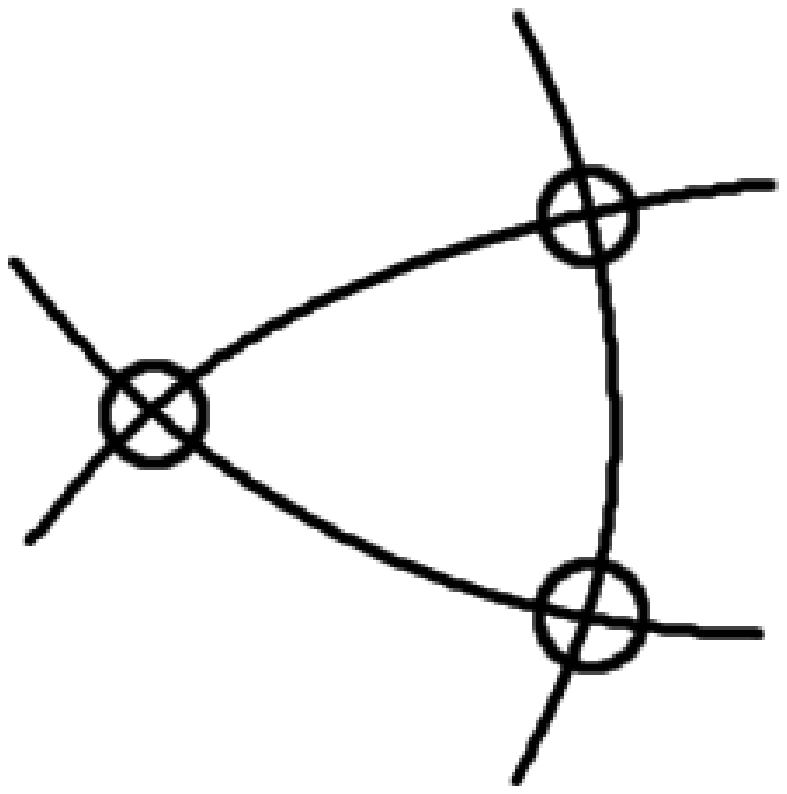}}
\raisebox{-0.5in}{\includegraphics[width=1in,height=1in]{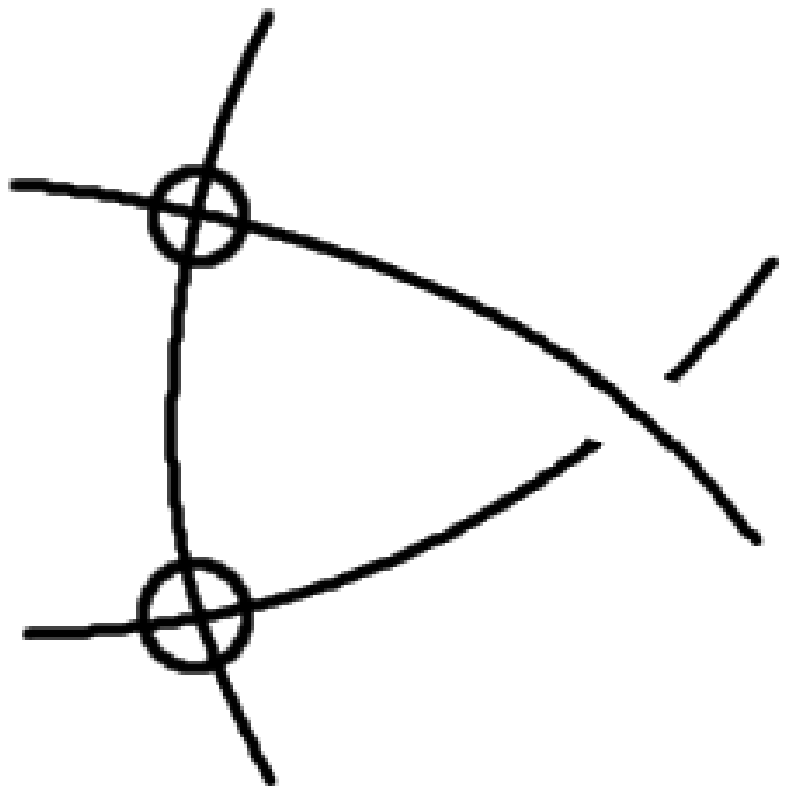}}
{\longleftrightarrow \atop v}
\raisebox{-0.5in}{\includegraphics[width=1in,height=1in]{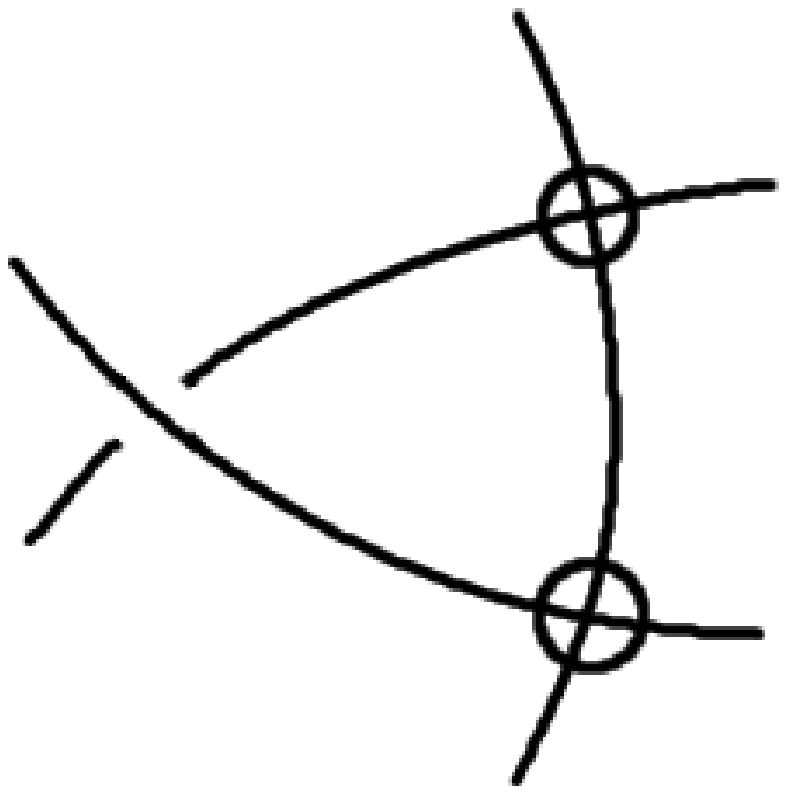}}
$$
\caption{Virtual Moves}
\end{figure}

Goussarov, Polyak and Viro [2] observe that there 
are two potential moves on virtual knot diagrams which resemble 
Reidemeister moves that are not 
allowed -- these ``forbidden moves'' depicted in figure 4 alter the Gauss 
diagram, unlike the other virtual moves. Worse yet, if these two moves 
are allowed, together they allow us to unknot any knot, rendering the theory 
trivial. For this reason, these moves are called ``forbidden''. On Gauss 
diagrams, the forbidden move $F_t$ moves an arrowtail of either sign
past an adjacent arrowtail with either sign without conditions 
on the relative positions of the heads of these arrows, and the other 
forbidden move $F_h$ moves an arrowhead past an arrowhead similarly.

\begin{figure}[h!]
$$ 
\raisebox{-0.5in}{\includegraphics[width=1in,height=1in]{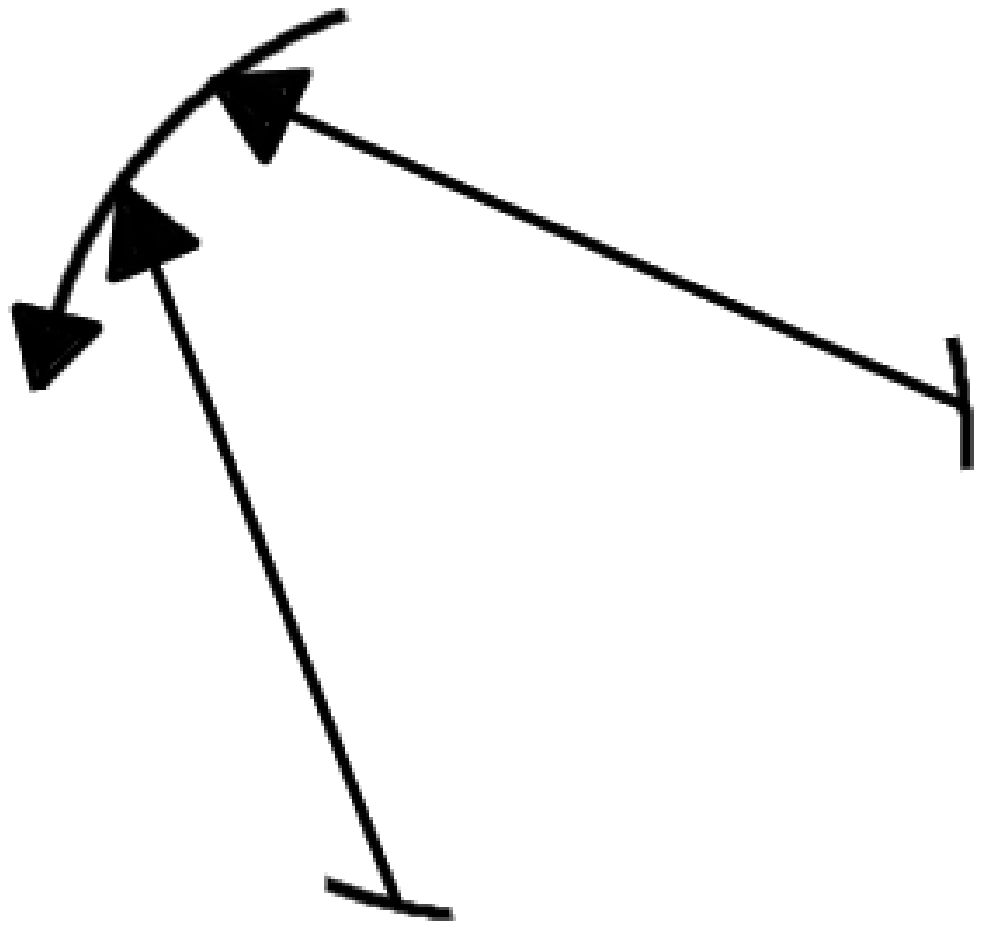}} 
{\longleftrightarrow \atop F_h}
\raisebox{-0.5in}{\includegraphics[width=1in,height=1in]{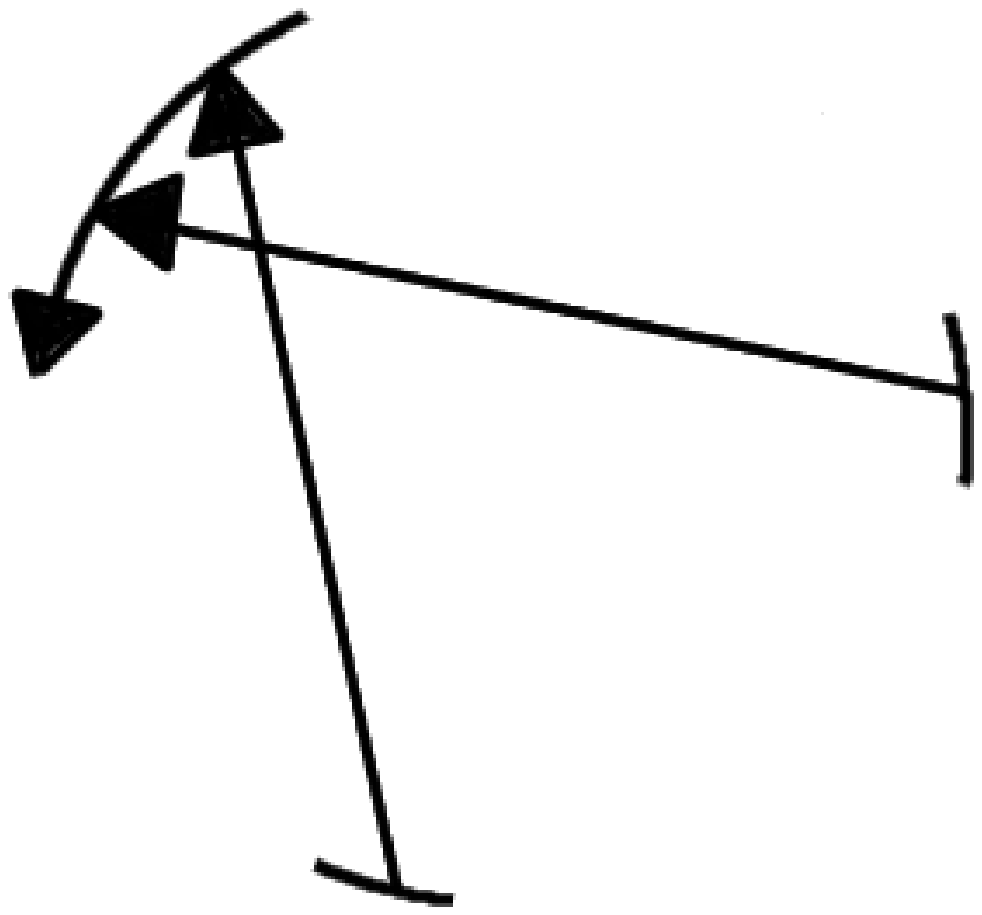}}
\raisebox{-0.5in}{\includegraphics[width=1in,height=1in]{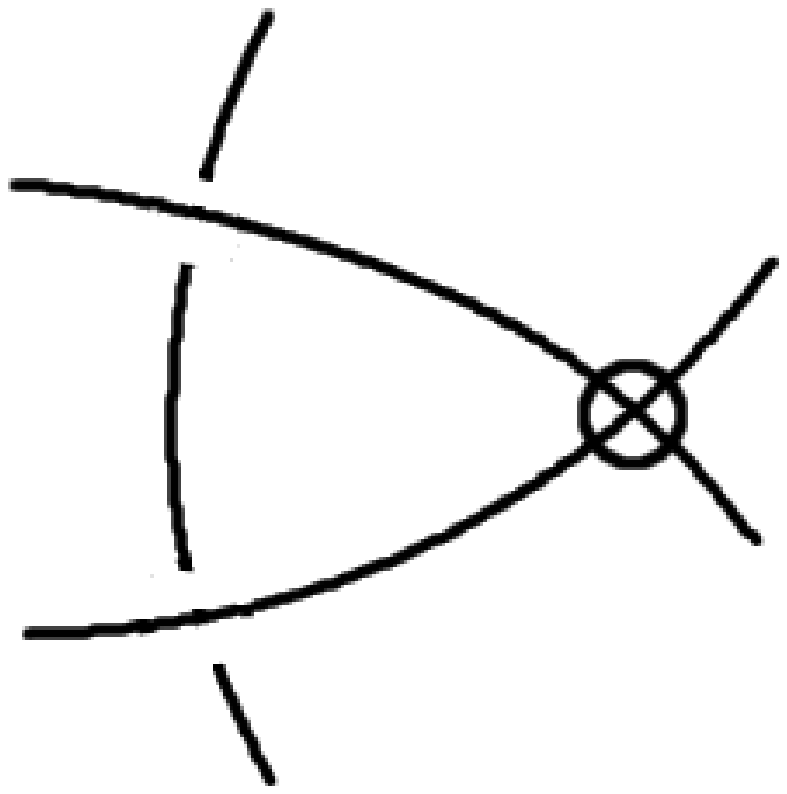}}
{\longleftrightarrow \atop F_h}
\raisebox{-0.5in}{\includegraphics[width=1in,height=1in]{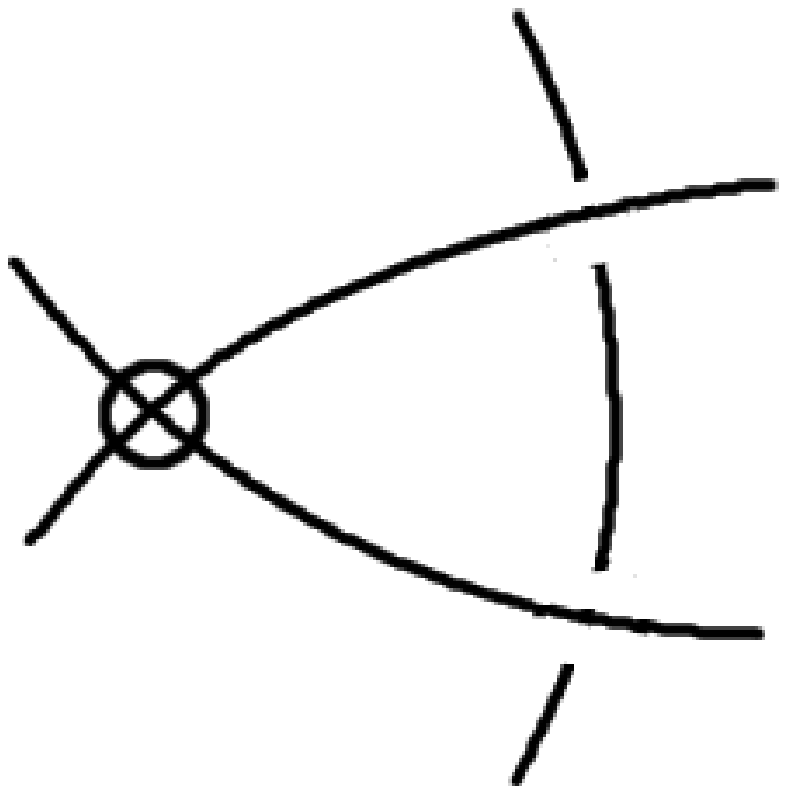}}
$$
$$
\raisebox{-0.5in}{\includegraphics[width=1in,height=1in]{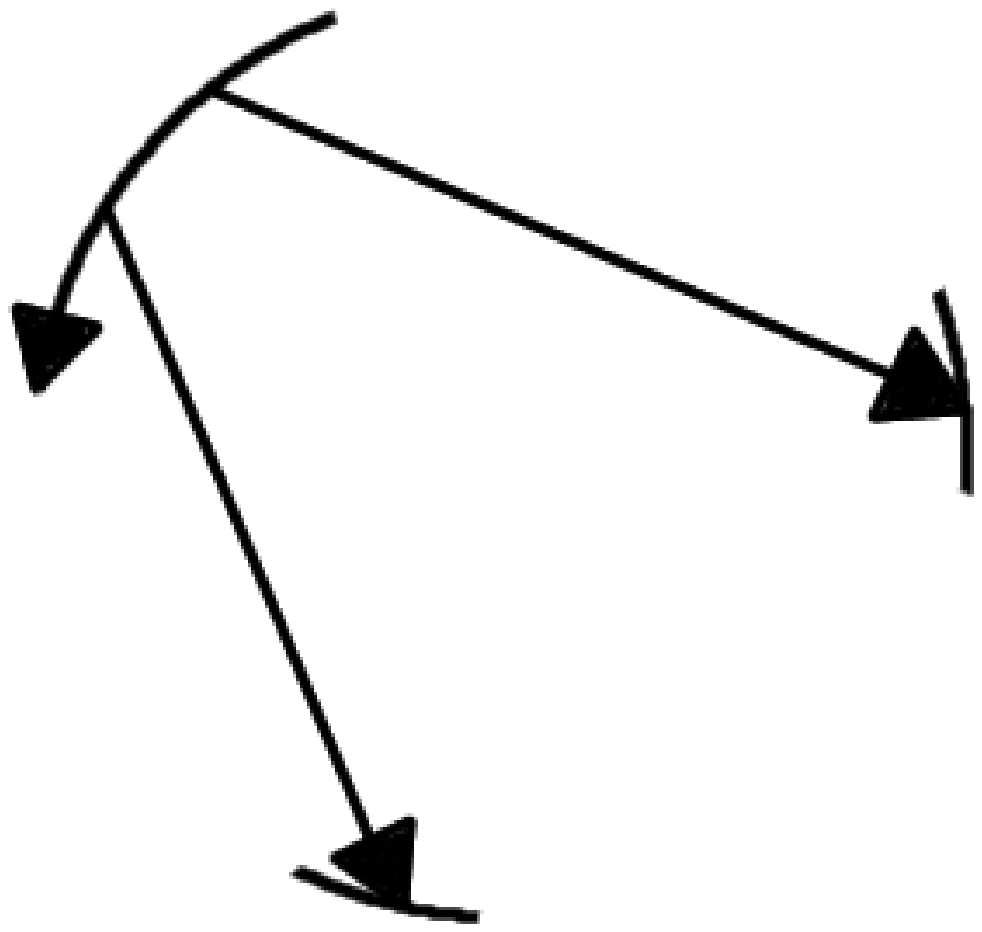}} 
{\longleftrightarrow \atop F_t}
\raisebox{-0.5in}{\includegraphics[width=1in,height=1in]{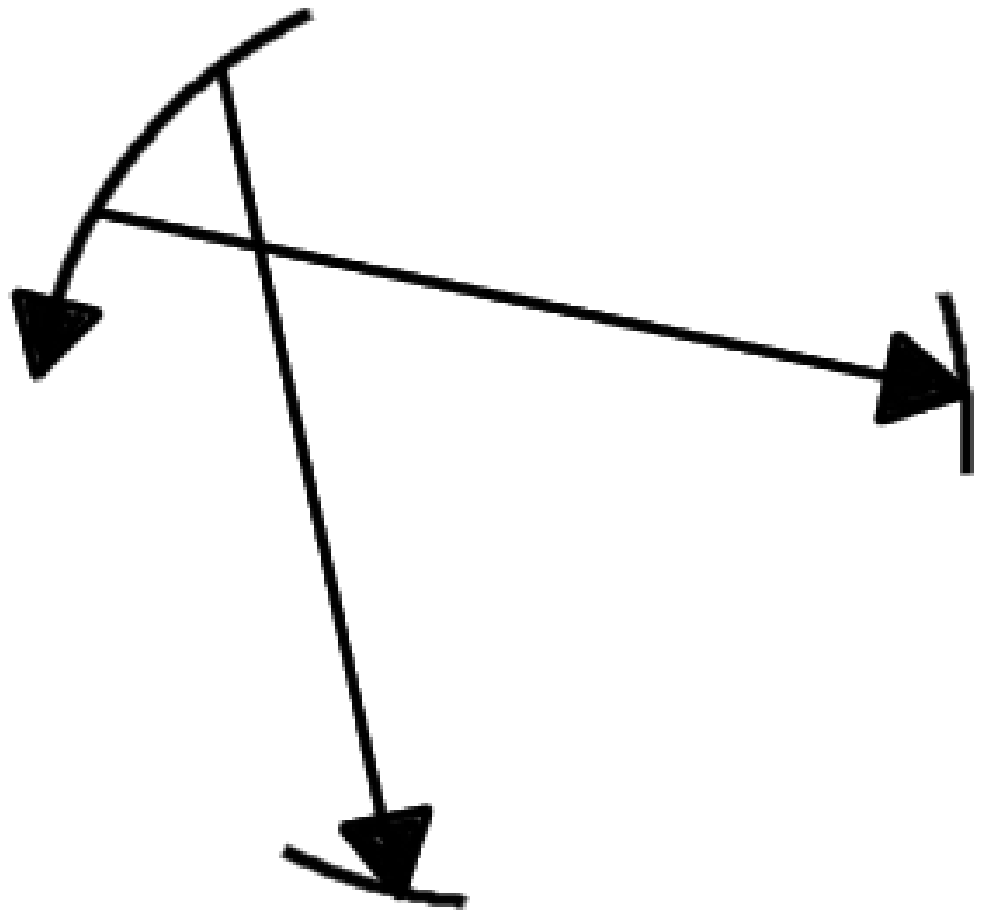}} 
\raisebox{-0.5in}{\includegraphics[width=1in,height=1in]{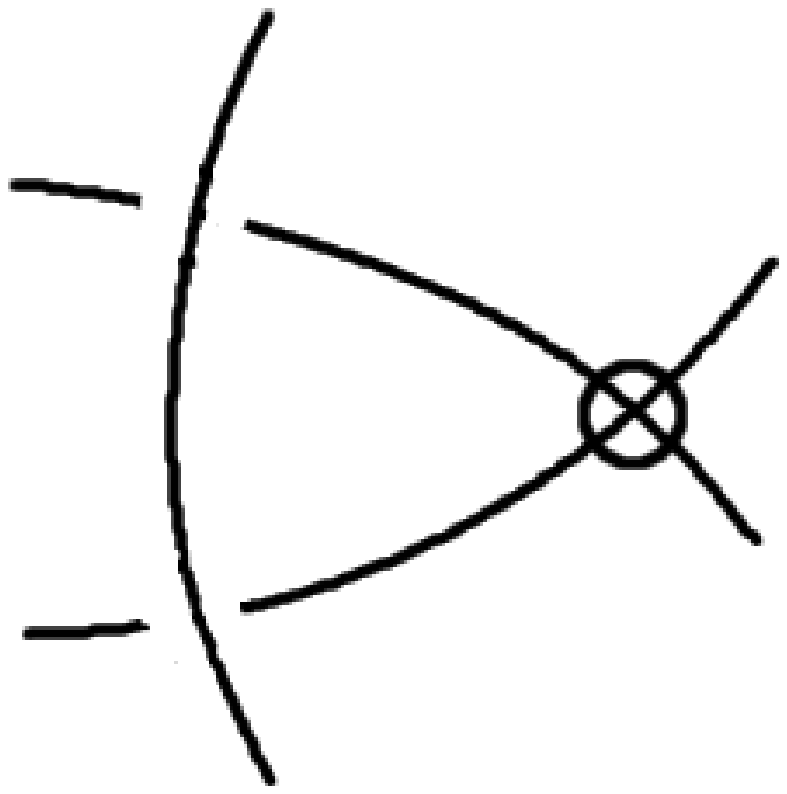}}
{\longleftrightarrow \atop F_t}
\raisebox{-0.5in}{\includegraphics[width=1in,height=1in]{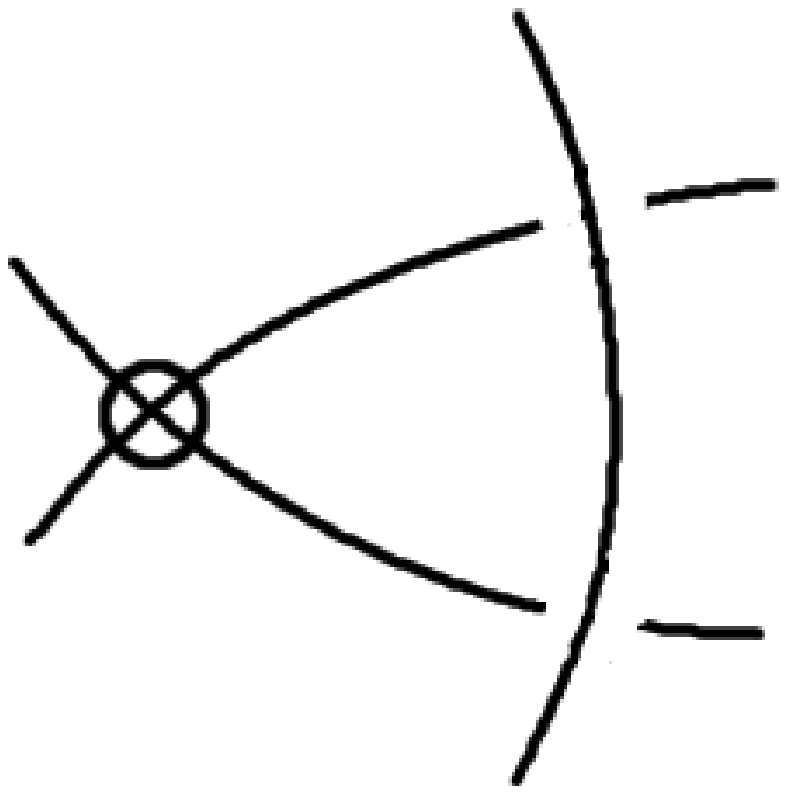}}
$$
\caption{The Forbidden moves $F_h$ and $F_t$.}
\end{figure}

The fact that allowing the forbidden moves would render virtual (and hence 
classical) knot theory trivial by making every knot unknotted is proven 
in [2] in terms of $n$-variations. In this paper we present a 
short combinatorial proof in terms of Gauss diagrams. The author has 
subsequently learned that Taizo Kanenobu [3] has a different
combinatorial proof of this result using virtual braid moves.

\medskip
\textbf{Theorem 1.}\textit{Any Gauss diagram can be changed into
any other Gauss diagram by a sequence of moves of types I, II, III, $F_h$ and
$F_t.$}

\medskip
\textbf{Proof.} The forbidden move $F_h$ allows us to move an 
arrowhead with either sign past an adjacent arrowhead with 
either sign without conditions on the tails of these respective arrows, 
and the move $F_t$ lets us do the same with arrowtails. If we could move 
an arrowhead of either sign past an arrowtail of either sign in the same 
manner, we could simply rearrange the arrows in a given diagram at will. 

The sequences of moves in figures 5 and 6 show how to move an arrowhead 
past an arrowtail of the same sign (move $F_s$; see Figure 4) or past an 
arrowtail of the opposite sign (move $F_o$; see Figure 5) using ordinary
Reidemeister moves and both forbidden moves.

\begin{figure}[h!]
$$ 
\raisebox{-0.5in}{\includegraphics[width=1in,height=1in]{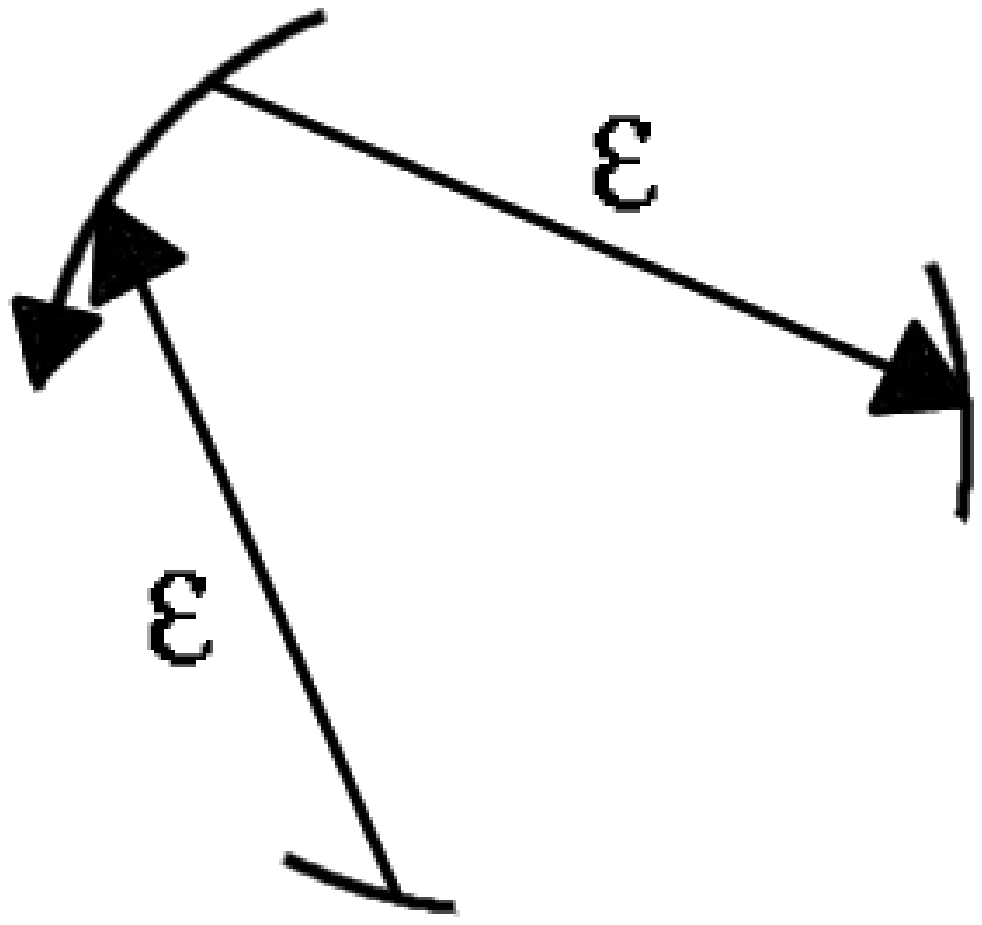}} 
{\longleftrightarrow \atop \mathrm{II}} 
\raisebox{-0.5in}{\includegraphics[width=1in,height=1in]{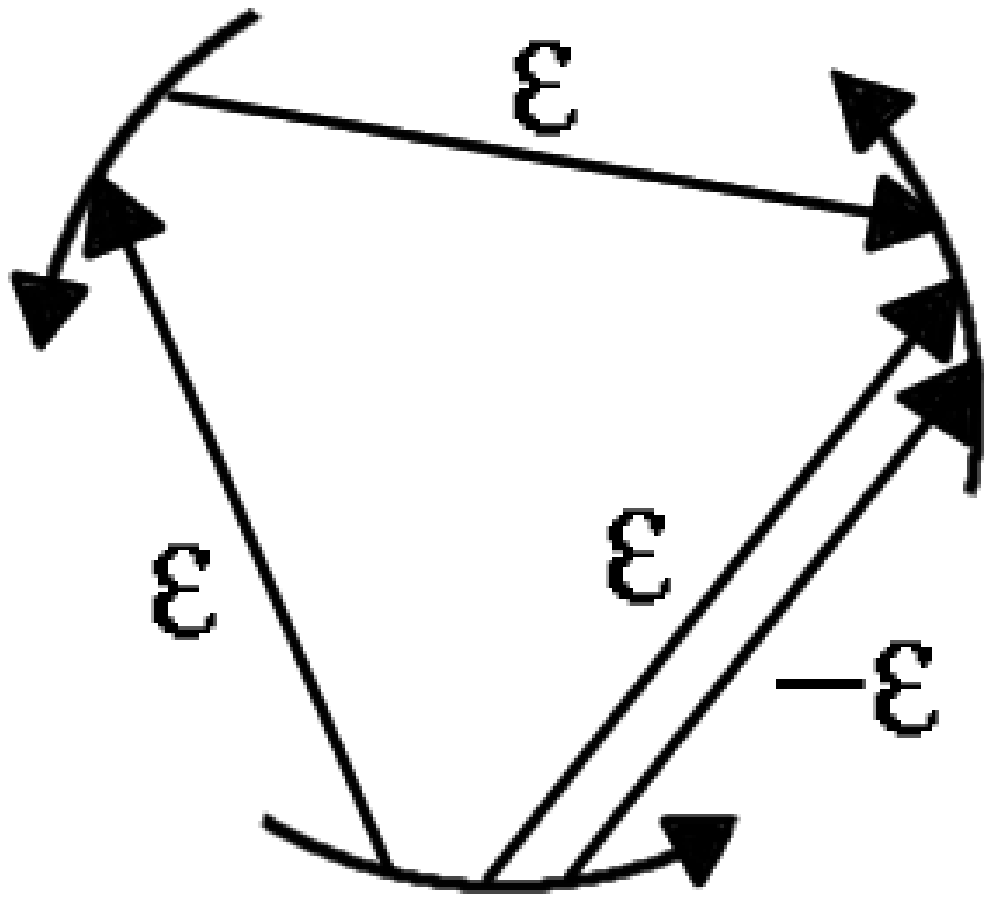}}
{\longleftrightarrow \atop F_t}
\raisebox{-0.5in}{\includegraphics[width=1in,height=1in]{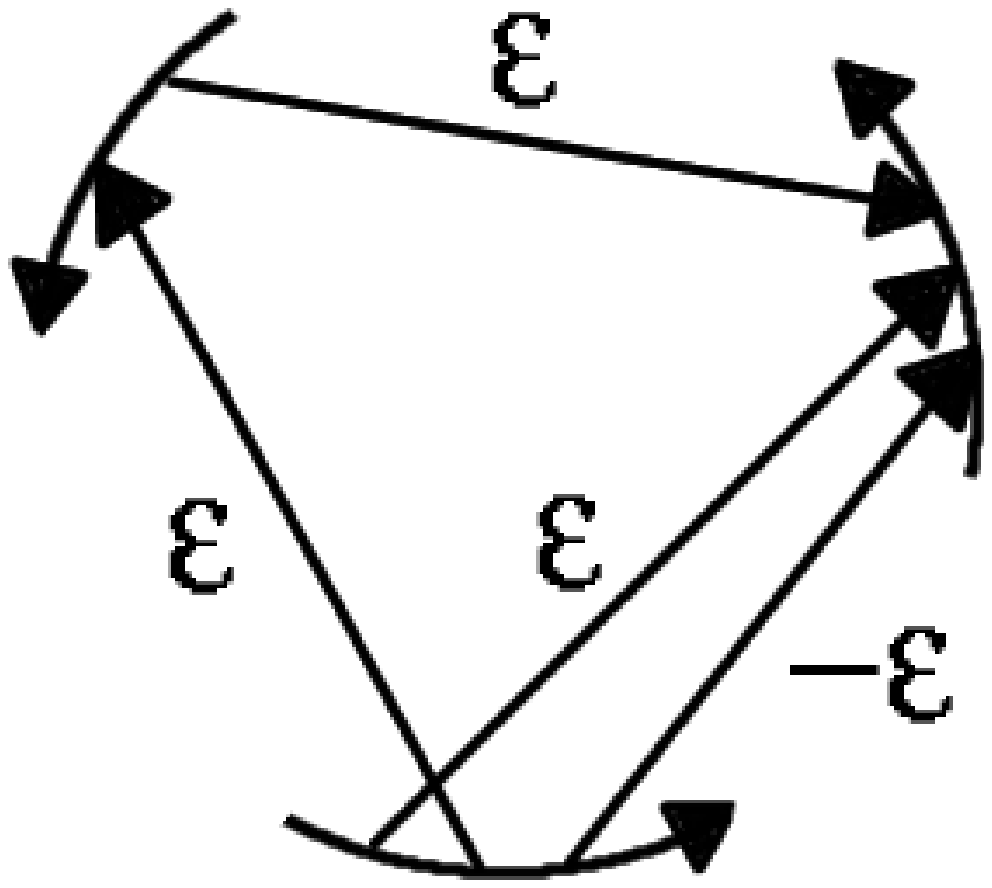}}
$$
$$
{\longleftrightarrow \atop F_h}
\raisebox{-0.5in}{\includegraphics[width=1in,height=1in]{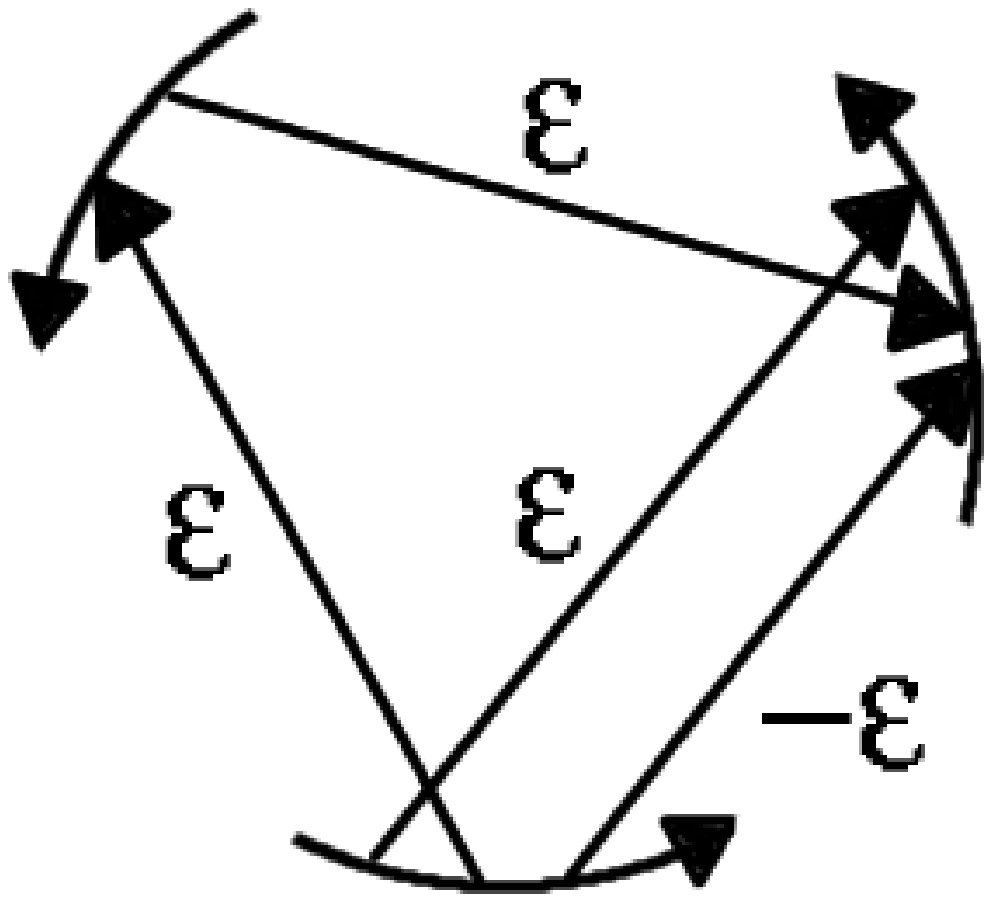}}
{\longleftrightarrow \atop \mathrm{III}}
\raisebox{-0.5in}{\includegraphics[width=1in,height=1in]{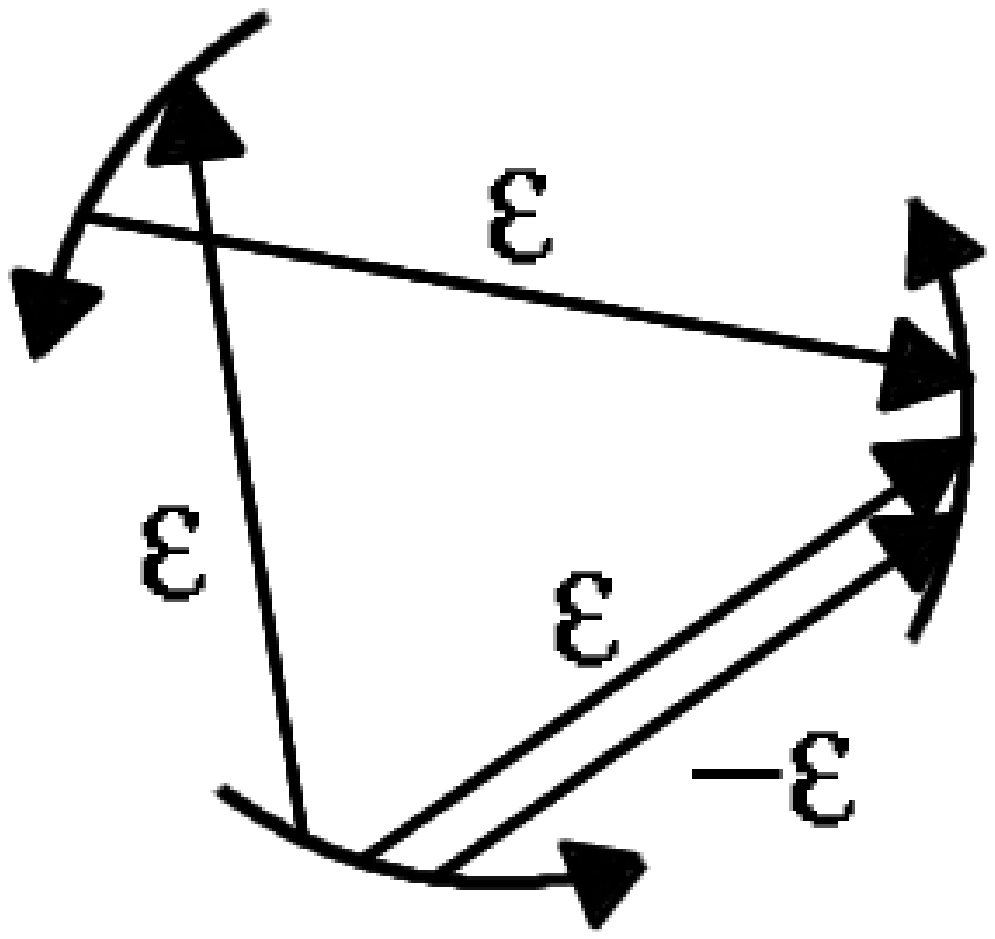}}
{\longleftrightarrow \atop \mathrm{II}}
\raisebox{-0.5in}{\includegraphics[width=1in,height=1in]{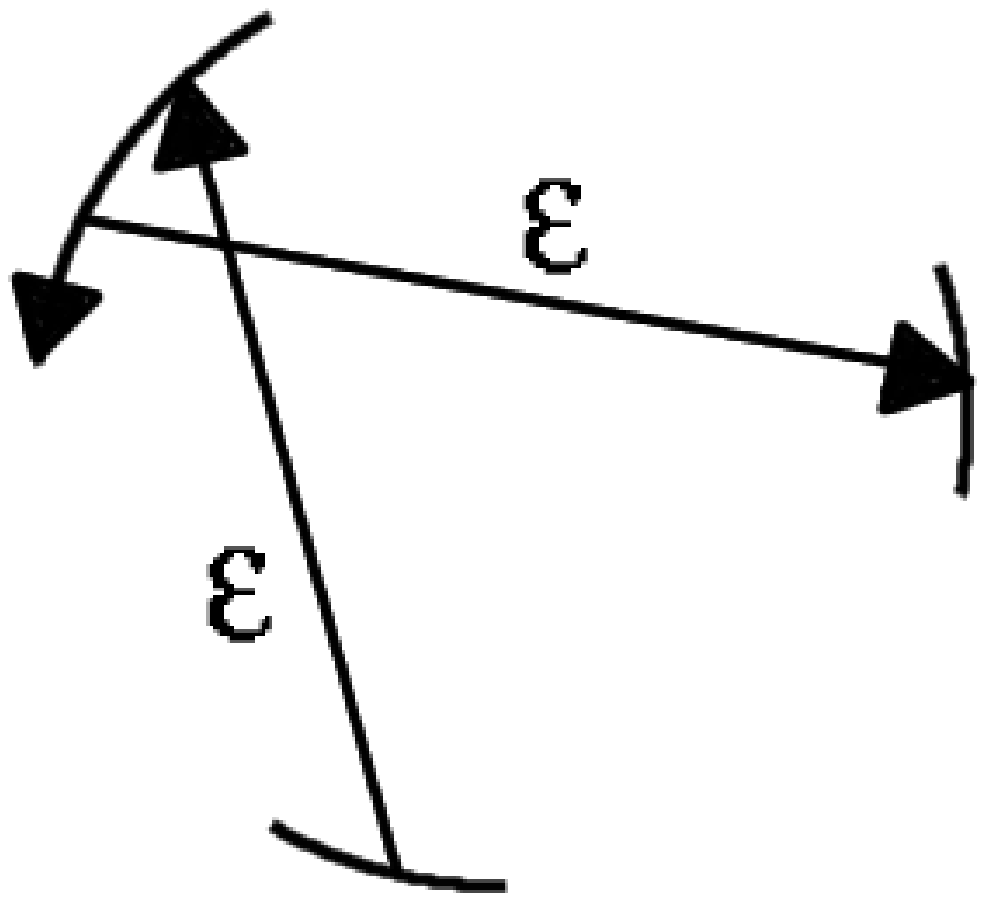}}
$$
\caption{Move sequence $F_s$.}
\end{figure}

\begin{figure}[h!]
$$ 
\raisebox{-0.5in}{\includegraphics[width=1in,height=1in]{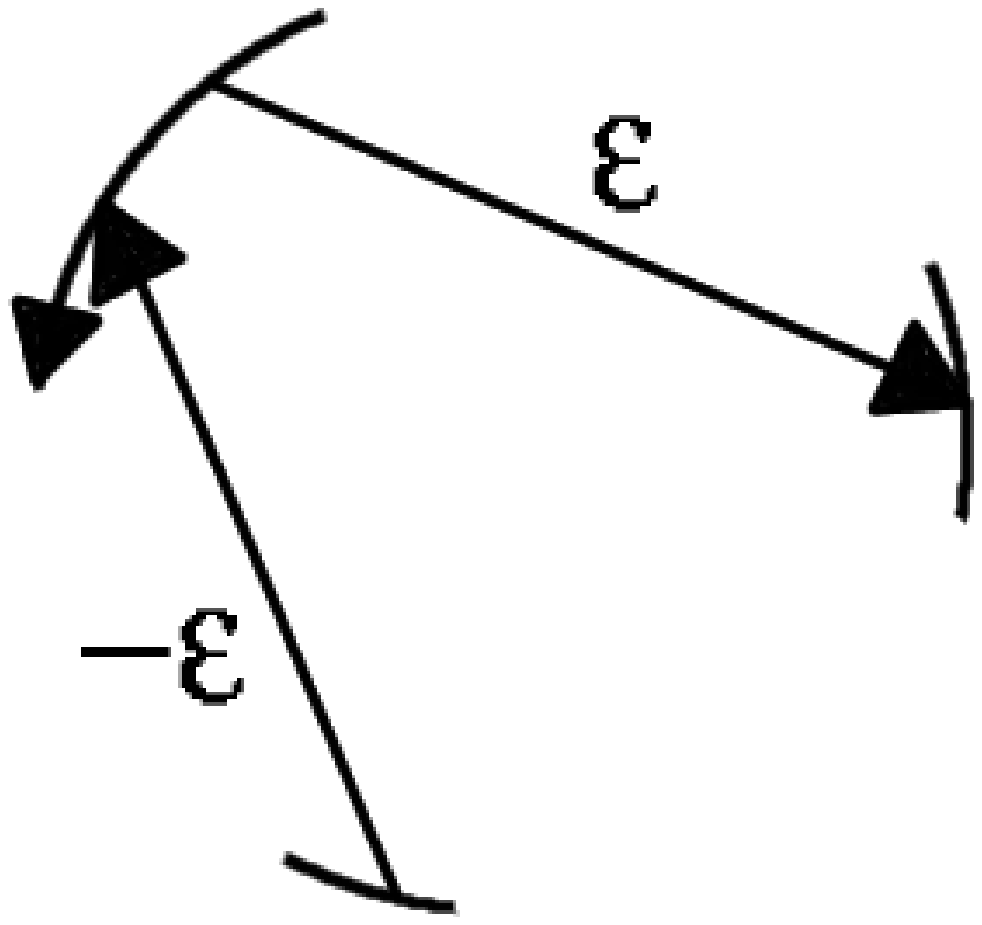}} 
{\longleftrightarrow \atop \mathrm{II}} 
\raisebox{-0.5in}{\includegraphics[width=1in,height=1in]{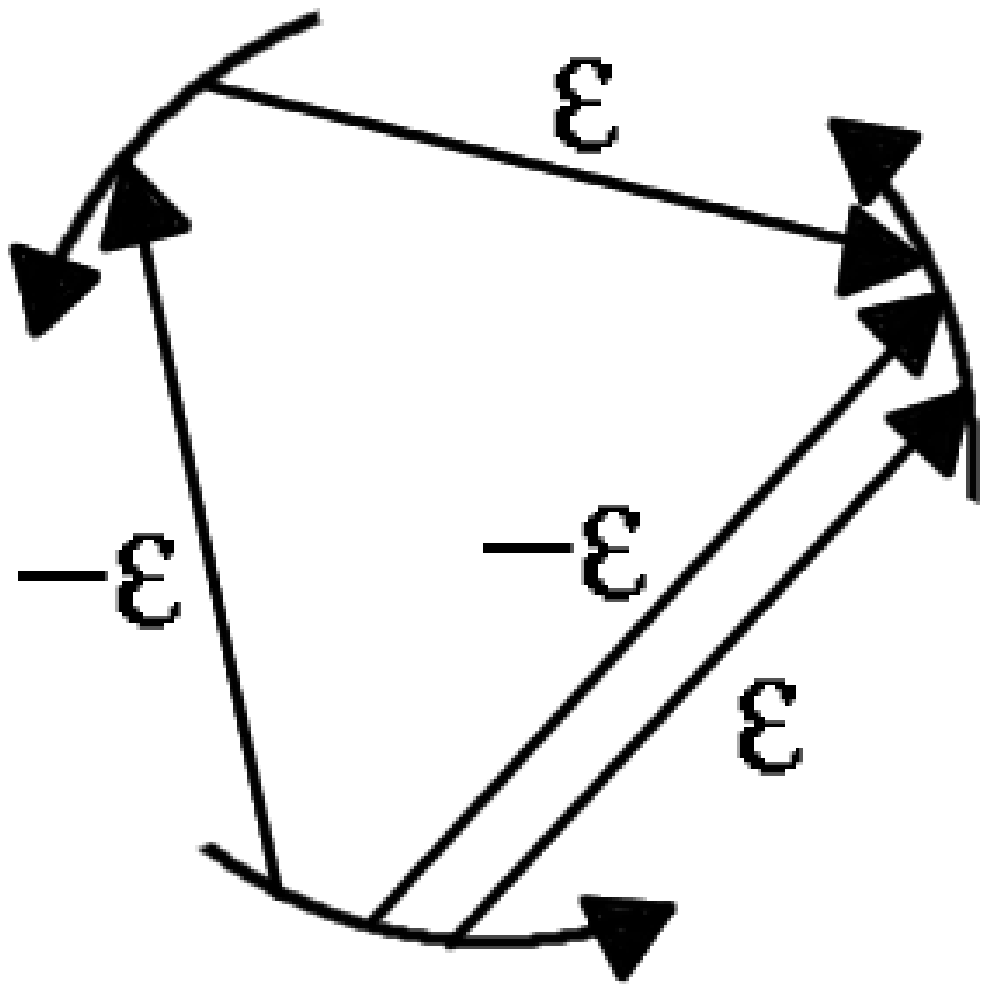}}
{\longleftrightarrow \atop F_h}
\raisebox{-0.5in}{\includegraphics[width=1in,height=1in]{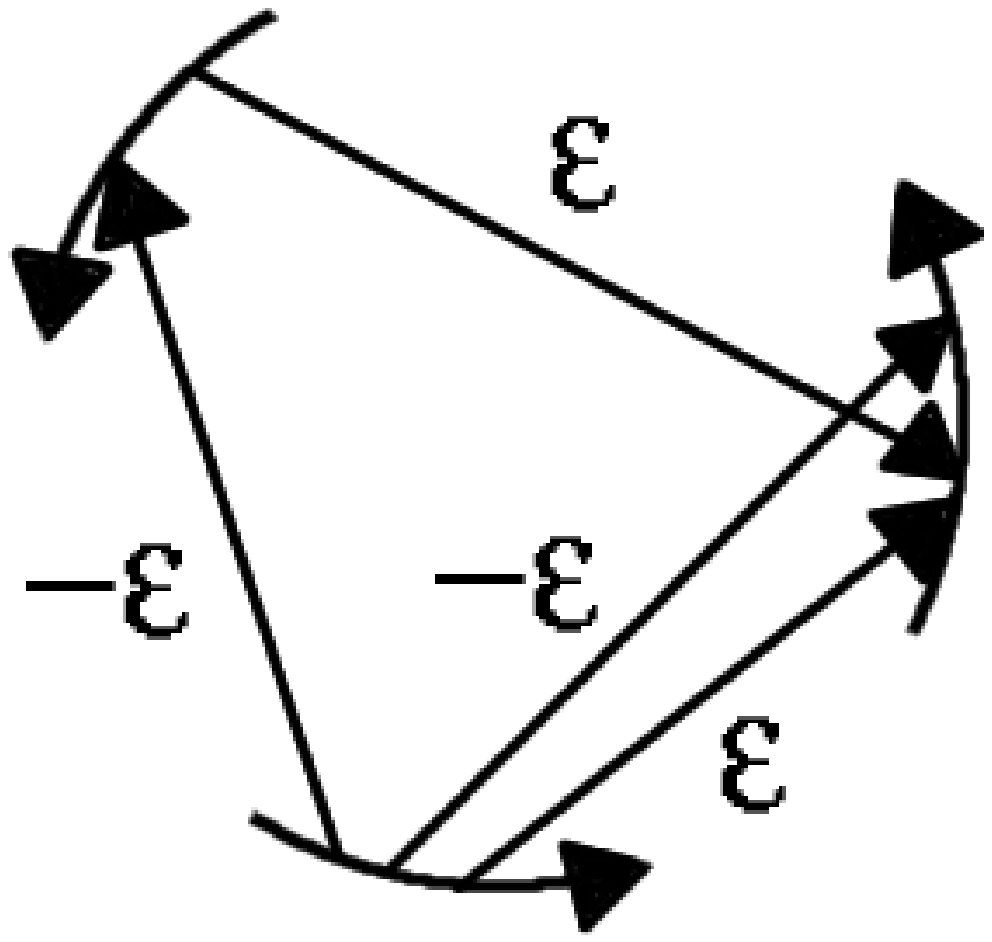}}
$$
$$
{\longleftrightarrow \atop \mathrm{III}}
\raisebox{-0.5in}{\includegraphics[width=1in,height=1in]{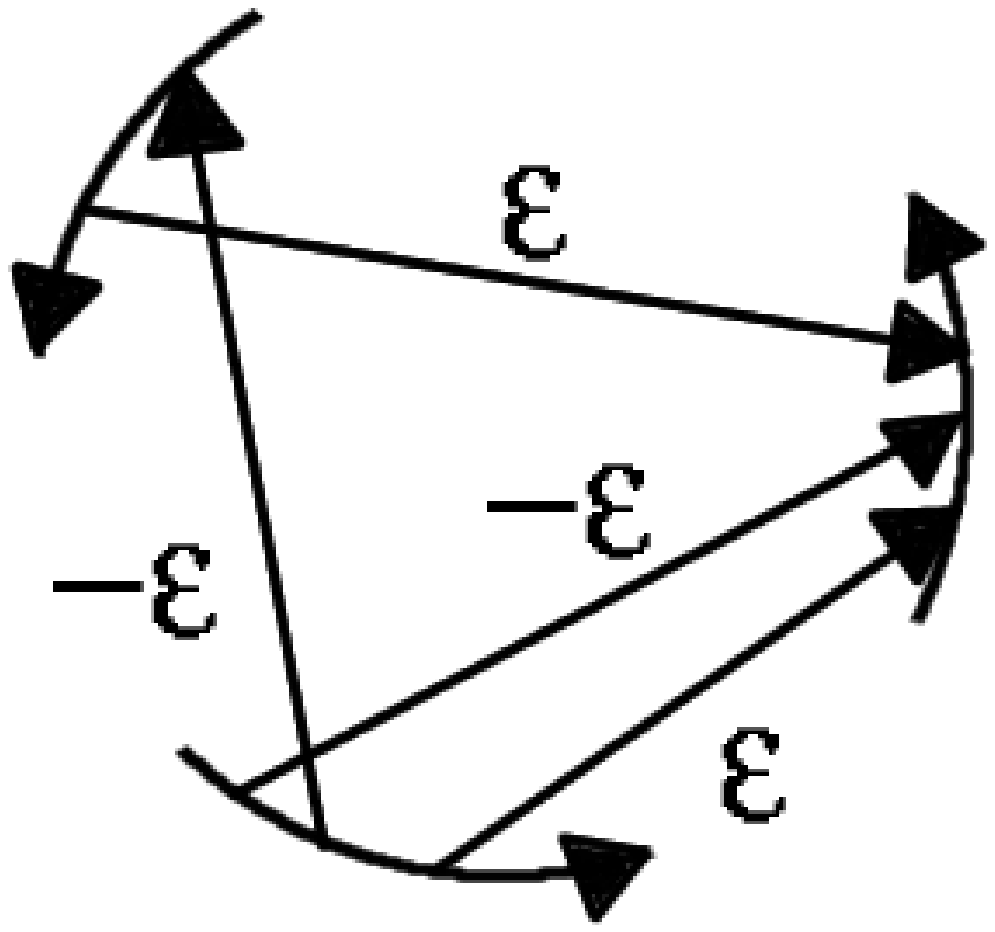}}
{\longleftrightarrow \atop F_t}
\raisebox{-0.5in}{\includegraphics[width=1in,height=1in]{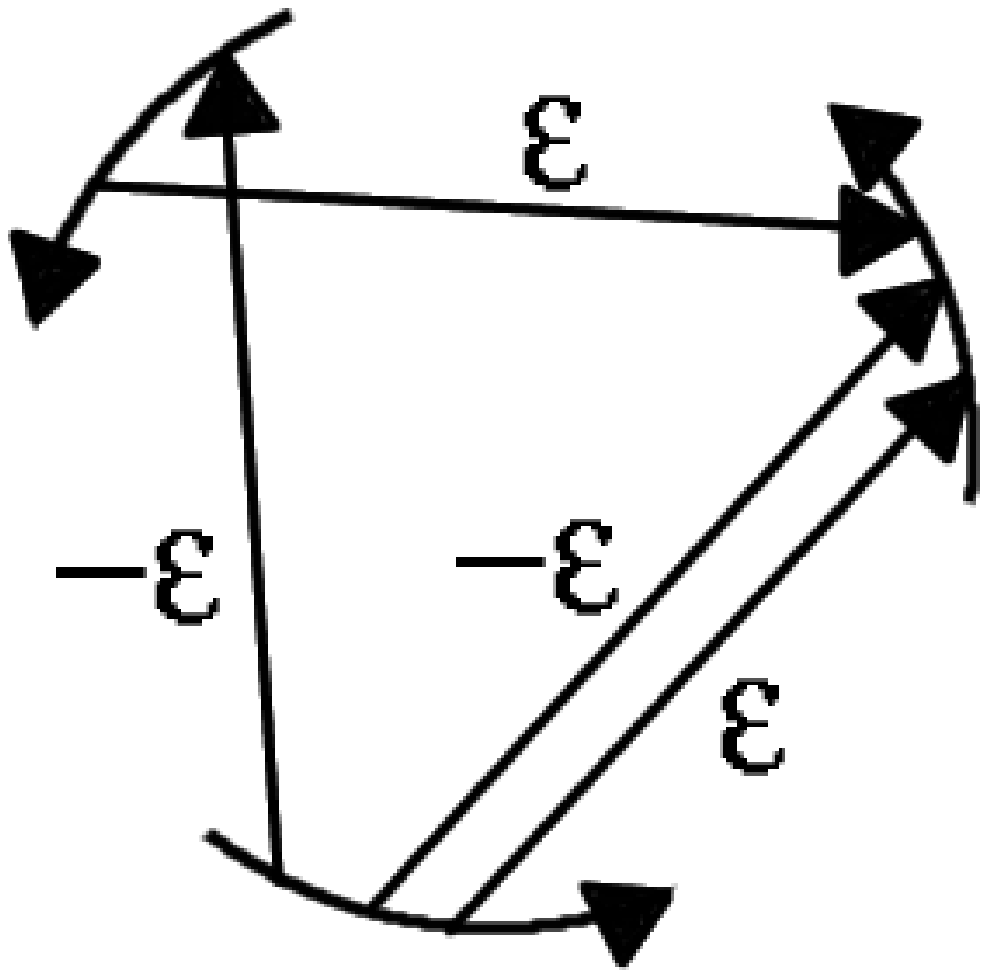}}
{\longleftrightarrow \atop \mathrm{II}}
\raisebox{-0.5in}{\includegraphics[width=1in,height=1in]{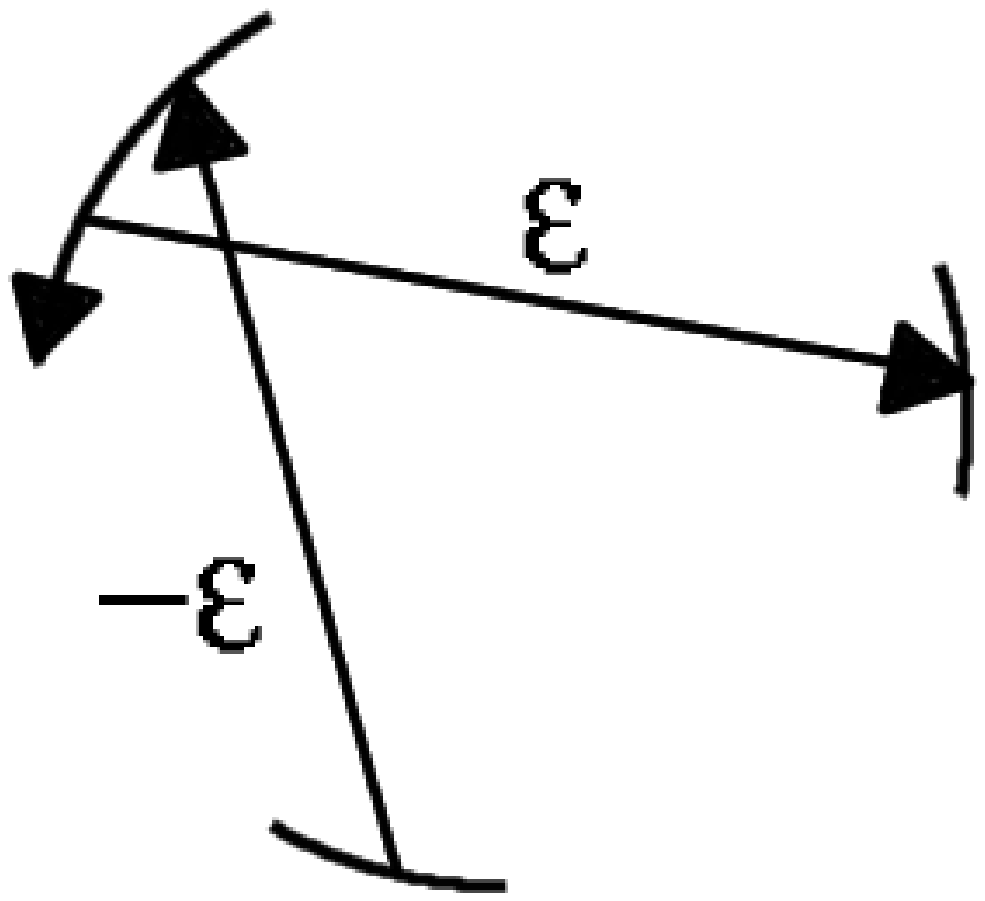}}
$$
\caption{Move sequence $F_o$.}
\end{figure}

Now, to change one Gauss diagram into another with the same
numbers of arrows of each sign, we simply use the forbidden moves and
moves sequences $F_o$ and $F_s$ to rearrange the arrows. If
we need more arrows of either sign, we can introduce them 
using type I moves and then move them into position with
moves the $F$ moves and sequences; if we have extra arrows, we 
can use the $F$ moves and sequences to move unwanted arrows into 
position to be removed by type I moves. In particular, any
virtual knot can be unknotted by this technique. $\qed $

\medskip
Note that the move sequences $F_s$ and $F_o$ each use both of the 
forbidden moves. If we define a new equivalence relation on Gauss
diagrams by allowing the usual Reidemeister moves and one forbidden move
but not the other, we arrive at the \textit{welded knots} of S. Kamada 
[4] or the \textit{weakly virtual knots} of S. Satoh [5]. Neither of the 
move sequences $F_s$ or $F_o$ can be constructed for 
welded knots since each move sequence requires the use of both
forbidden moves.

\bigskip
\noindent\textbf{References}
\medskip

{\small
\noindent [1] L. H. Kauffman, {\em Virtual Knot Theory}, 
Europ. J. Combinatorics 20 (1999) 663-690.

\noindent [2] M. Goussarov, M. Polyak and O.Viro, {\em 
Finite type invariants of classical and virtual knots},
preprint (http://xxx.lanl.gov/abs/math.GT/98100073).

\noindent [3] T. Kanenobu, {\em Forbidden moves unknot
a virtual knot}, J. Knot Theory Ramifications (to appear).

\noindent [4] S. Kamada, {\em Braid presentation of
virtual knots and welded knots}, preprint.

\noindent[5] S. Satoh, {\em Virtual knot presentation
of ribbon torus-knots}, J. Knot Theory Ramifications (to appear).
}

\end{document}